\documentclass[12pt]{amsart}

\usepackage{amsmath,amssymb,epic,lscape}

\newtheorem{theorem}{Theorem}[section]

\theoremstyle{definition}

\theoremstyle{remark}

\numberwithin{equation}{section}
\providecommand{\bysame}{\leavevmode\hbox to3em{\hrulefill}\thinspace}

\def\DJ{{\hbox{D\kern-.8em\raise.15ex\hbox{--}\kern.35em}}}
\def\DJo{$\;$\kern-.4em
    \hbox{D\kern-.8em\raise.15ex\hbox{--}\kern.35em okovi\'c}}

\def\NSERC{Supported in part by an NSERC Discovery Grant.}

\renewcommand{\subjclassname}{\textup{2000} Mathematics Subject
Classification }

\begin{document}

\title[A new Yang number]
{A new Yang number and consequences}

\author[D.\v{Z}. \DJ okovi\'{c}]
{Dragomir \v{Z}. \DJ okovi\'{c}} 
\address{Department of Pure Mathematics, University of Waterloo,
Waterloo, Ontario, N2L 3G1, Canada}

\email{djokovic@uwaterloo.ca}

\thanks{\NSERC}


\date{}

\begin{abstract}
Base sequences $BS(m,n)$ are quadruples $(A;B;C;D)$ of 
$\{\pm1\}$-sequences, $A$ and $B$ of length $m$ and $C$ and $D$
of length $n$, the sum of whose non-periodic auto-correlation
functions is zero. Base sequences and some special subclasses of 
$BS(n+1,n)$ known as normal and near-normal sequences, 
$NS(n)$ and $NN(n)$, as well as T-sequences and orthogonal
designs play a prominent role in modern constructions of
Hadamard matrices. In our previous papers \cite{DZ2,DZ3} we have 
classified the near-normal sequences $NN(s)$ for all even integers 
$s\le32$ (they do not exist for odd $s>1$). 
We now extend the classification to the case $s=34$.
Moreover we construct the first example of near-normal sequences 
$NN(36)$. Consequently, we construct for the first time 
T-sequences of length 73. For all smaller lengths, T-sequences
were already known. Another consequence is that 73 is a 
Yang number, and a few important consequences of this fact 
are given.
\end{abstract}

\maketitle
\subjclassname{ 05B20, 05B30 }
\vskip5mm

\section{Preliminaries}

A sequence $A=a_1,a_2,\ldots,a_m$, of length $m$, is {\em binary} 
respectively {\em ternary} if $a_i\in\{\pm1\}$ respectively $a_i\in\{0,\pm1\}$.
We identify $A$ with the polynomial
$A(z)=a_1+a_2z+\cdots+a_mz^{m-1}$. The {\em norm} of a Laurent 
polynomial $f(z)$ is defined as $N(f)=f(z)f(z^{-1})$.
A quadruple $(A;B;C;D)$ of binary sequences, with $A$ and $B$
of length $m$ and $C$ and $D$ of length $n$, are {\em base sequences} if 
$$ N(A)+N(B)+N(C)+N(D)=2(m+n). $$
The set of such base sequences is denoted by $BS(m,n)$.
The {\em Base Sequence Conjecture} (BSC) asserts that 
$BS(n+1,n)\ne\emptyset$ for all nonnegative integers $n$
(see \cite{DZ1}). It has been confirmed for $n\le35$. It is 
also well known that it holds when $n$ is a {\em Golay number},
i.e., a number of the form 
$2^a 10^b 26^c$ where $a,b,c$ are nonnegative integers.

The base sequences $(A;B;C;D)\in BS(n+1,n)$ are {\em normal}
respectively {\em near-normal} if $a_i=b_i$ respectively $a_i=(-1)^{i-1}b_i$ 
for all indexes $i$ in the interval $1\le i\le n$.
Let $NS(n)$ resp. $NN(n)$ be the subset of $BS(n+1,n)$
consisting of all normal resp. near-normal sequences.
It is known that $NS(n)\ne\emptyset$ when $n$ is a Golay 
number. For positive integers $n\le33$, $NS(n)=\emptyset$
exactly for 
$$ n\in \{6,14,17,21,22,23,24,27,28,30,31,33\} $$
(see \cite{DZ3}).
It is known that $NN(n)=\emptyset$ when $n$ is odd.
{\em Yang Conjecture} (YC) asserts that $NN(n)\ne\emptyset$ 
when $n$ is even. This has been confirmed for $n\le34$
(see \cite{DZ3}).

A quadruple of ternary sequences $(A;B;C;D)$, all four
of length $n$, are {\em $T$-sequences} if
$$ N(A)+N(B)+N(C)+N(D)=n $$
and for each index $i$, $1\le i\le n$, exactly one of 
$a_i,b_i,c_i,d_i$ is nonzero.

Let $x_i$, $i=1,2,\ldots,u$, be independent commuting
variables. An $n\times n$ matrix $S$ all of whose entries
belong to $\{0,\pm x_1,\ldots,\pm x_u\}$ and such that
$$ SS^T=(s_1x_1^2+\cdots+s_ux_u^2)I_n $$
is an {\em orthogonal design}. Here the superscript $T$
denotes transposition of matrices and $I_n$ is the
identity matrix. The symbols $s_1,\ldots,s_u$ are
positive integers. The set of all such orthogonal
designs is denoted by $OD(n;s_1,\ldots,s_u)$.

Yang \cite[Theorems 1 and 3]{CY} constructed two maps,
known as {\em Yang multiplications},
\begin{eqnarray*}
&& NS(s) \times BS(m,n) \to TS((2s+1)(m+n)), \\
&& NN(s) \times BS(m,n) \to TS((2s+1)(m+n)).
\end{eqnarray*}
(As pointed out in \cite{DZ3}, there are two misprints in 
his Theorem 1.) A {\em Yang number} is an odd integer $2s+1$ 
such that $NS(s)$ or $NN(s)$ is not empty.

\section{Main result and consequences}

We have carried out an exhaustive search for $NS(n)$ when 
$n=34,35$ and found that both are empty.

In our paper \cite{DZ2} we have introduced an equivalence
relation, NN-equivalence, and a canonical 
form for near-normal sequences and we have enumerated
the equivalence classes in $NN(n)$ for even integers $n\le 30$.
For the classification in the case $n=32$ see \cite{DZ3}.
We have now completed the classification in the case $n=34$.
The representatives $(A;B;C;D)$ of the NN-equivalence classes
of $NN(34)$ are given in Table 1. The sums $a,b,c,d$ of the
sequences $A,B,C,D$ are also recorded.
The sequences are given in encoded form for the sake
of compactness. The encoding scheme is explained in our
papers \cite{DZ1,DZ2}. This table also contains
the unique known example (up to NN-equivalence) of $NN(36)$,
which we have constructed recently.

\begin{theorem} \label{glavna}
$NN(36)\ne\emptyset.$
\end{theorem}
\begin{proof} Direct verification. \end{proof}

\begin{center}
\begin{tabular}{|r|l|l|r|}
\multicolumn{4}{c}{\bf Table 1: Near-normal sequences $NN(n)$} \\
\multicolumn{4}{c}{} \\ \hline 
\multicolumn{1}{|c|}{} & \multicolumn{1}{c|}{$A$ \& $B$} & 
\multicolumn{1}{c|}{$C$ \& $D$} & \multicolumn{1}{c|}{$a,b,c,d$} \\ \hline
\multicolumn{4}{c}{$n=34$} \\ \hline
1 & 076417646512321462 & 16738541372344337 & $7,7,-2,6$  \\	
2 & 076535878535141762 & 17677852174231455 & $-5,7,0,8$ \\	
3 & 076782178767646231 & 17621532262576812 & $-5,3,10,-2$ \\
4 & 058214353712141461 & 11868756376664254 & $11,3,-2,2$ \\	
5 & 053765656464871261 & 17765746348615187 & $1,1,-6,10$ \\	
\hline 
\multicolumn{4}{c}{$n=36$} \\ \hline
1 & 0764841234846532153 & 165154775335162126 & $3,-3,8,8$ \\	
\hline
\multicolumn{4}{c}{} \\
\end{tabular}
\end{center}

For the convenience of the reader we decode the above example of 
$NN(36)$ and write it as a quadruple $(A;B;C;D)$ of binary sequences. 
In fact we write $+$ for $+1$ and $-$ for $-1$.

\begin{eqnarray*}
A&=& +-++-+++-+-++--++--+ \\
&\ & ++++----+++-----+; \\
B&=& +++---+-----++--++-- \\
&\ & +-+--+-++-++-+-+-; \\
C&=& ++-+-+------++++++-+ \\
&\ & ++-+++++---+++-+; \\
D&=& +++++-+++--+++-+-+-- \\
&\ & +--+-++-+++-+--+. \\
\end{eqnarray*}

Let us list some consequences of the above results.

(i) Since $NN(n)\subseteq BS(n+1,n)$, we can update the status
of the BSC and YC: $BS(n+1,n)\ne\emptyset$ for $n\le36$, and
$NN(n)\ne\emptyset$ for even $n\le36$.

(ii) An odd (positive) integer $n\le 73$ is a Yang number 
if and only if $n\ne 35,43,47,55,63,67,71.$ In order to rule out the
integer $71$, we need to use the fact that $NS(35)=\emptyset$ 
mentioned above.

(iii) Our main result implies that there exist T-sequences $TS(73)$ 
(see \cite[Lemma 5.21]{SY}). Consequently, by \cite[Theorem 3.6]{SY}, 
there exists an orthogonal design $OD(4t;t,t,t,t)$ for $t=73$. 
Neither T-sequences nor T-matrices for $n=73$ were known previously, 
see the Remarks V 2.51 and V 8.47 in \cite{HCD}. 
In spite of the claim made in Remark V 2.119.5 that an OD with the 
above parameters is known, we believe that this is not the case as we
could not find such a result anywhere in the literature and our request 
for a reference failed.

(iv) One can plug into our OD any Williamson-type matrices 
of some order $n$ to obtain Hadamard matrices of order $4\cdot73\cdot n$. 
Since Williamson-type matrices are known for infinitely many odd orders 
$n$, we obtain infinitely many Hadamard matrices of order an odd 
multiple of 4. In particular, there exist Williamson-type matrices of 
order $n=61,81,83,101$ (see \cite[Table A.1]{SY}). Hence, there exist 
Hadamard matrices of order $4\cdot73\cdot n$ for the same values of $n$.
These orders fall in the range covered by Table V 1.53 of \cite{HCD}, 
where it is indicated that no Hadamard matrices of these orders are
known. However, we have recently discovered over 100 errors in this
table. It turns out that Hadamard matrices of these four orders are
in fact known. This will be discussed in more details elsewhere.

(v) For $s=2,4,\ldots,34$, the number of NN-equivalence classes in 
$NN(s)$ is $1,2,2,3,8,14,11,24,20,18,32,12,3,20,9,8,5.$

Our computer program performs an exhaustive search for near-normal sequences
$NN(s)$ for fixed $s$. The search is divided into 12 cases which can be
run separately. For $s=36$, we started to run them at different times 
according to the availability of machines. Six of the cases have
completed within a month or two without finding any near-normal
sequences. Only one of the remaining cases produced (after about 40 days)
the example given in the above table.

It is a pleasure to thank an anonymous referee for his valuable
comments.

\end{document}